 \documentclass{article}

\usepackage{amssymb}
\usepackage{amsmath}
\usepackage{amsfonts}
\usepackage{pst-all}
\usepackage{pst-poly}
\usepackage{graphics}
\usepackage{graphicx}
\usepackage{ntheorem}
\usepackage{tabularx}

\newcounter{num}

\newcommand{\rr}{\mathbb{R}}
\newcommand{\zz}{\mathbb{Z}}
\newcommand{\cc}{\mathbb{C}}

\newcommand{\dd}{\mathbb{D}}
\newcommand{\sss}{\mathbb{S}}

\newcommand{\ttt}{\mathbb{T}}

\newtheorem{teorema}{Theorem}[section]
\newtheorem{definicion}[teorema]{Definition}
\newtheorem{lema}[teorema]{Lemma}
\newtheorem{proposicion}[teorema]{Proposition}

\theoremnumbering{Alph}
\newtheorem{Teorema}{Theorem}

\newcommand{\demostracion}[1]{
		\noindent\textbf{Proof.}
		#1 \par
	       \hfill $\Box$
}

\newcommand{\sdemostracion}[1]{\vspace{1mm}
		\noindent\textbf{Proof}
		#1 \par
	       \hfill $\Box$
}

\begin{document}

\title{Existence of periodic orbits for geodesible vector fields on closed 3-manifolds}

\author{Ana Rechtman}

\maketitle
\begin{abstract}
In this paper we deal with the existence of periodic orbits of
geodesible vector fields on closed 3-manifolds. A vector field is
geodesible if there exists a Riemannian metric on the ambient manifold
making its orbits geodesics. In particular, Reeb vector fields and
vector fields that admit a global section are geodesible. We will
classify the closed 3-manifolds that admit aperiodic volume preserving $C^\omega$ geodesible vector
fields, and prove the existence of periodic orbits for
$C^\omega$ geodesible vector fields (not volume preserving), when the 3-manifold is not a
torus bundle over the circle. We will also prove the existence of
periodic orbits of $C^2$ geodesible vector fields in some closed
3-manifolds.
\end{abstract}

\section{Introduction}
\label{kristinahofer}

The already disproved Seifert's conjecture stated the existence of a periodic orbit for every non singular vector field on $\sss^3$. We can still talk about Seifert's conjecture on some families of non singular vector fields: for example {\it geodesible} or volume preserving.  From now on, if not otherwise stated, we will write flow for non singular flow, and vector field for non singular vector field. 

\begin{definicion}\label{defgeodesibles}
Let $X$ be a vector field on a closed manifold $M$, we will say that it is geodesible, or that the associated  flow is geodesible, if there exists a Riemannian metric $g$ on $M$ such that the orbits are geodesics. 
\end{definicion}

In this paper we are mainly interested in real analytic geodesible vector fields, in this case we will assume that the Riemannian metric is also real analytic. In general, when we talk about $C^r$ geodesible vector fields we will assume that there exists a $C^r$ Riemannian metric making its orbits geodesics. Let $X$ be a geodesible vector field and $g$ the Riemannian metric making its orbits geodesics. Then we have a 1-form $\alpha=\iota_Xg$, such that modulo reparameterization of $X$
$$
\alpha(X)=1 \qquad \mbox{and} \qquad \iota_Xd\alpha=0.
$$
Conversely, if $X$ is a vector field and $\alpha$ is a one form satisfying the above equations, $X$ is geodesible. We will prove the following two theorems

\begin{Teorema}\label{teoremaperiogeodanavol}
Let $X$ be a real analytic geodesible vector field on an oriented closed 3-manifold $M$, that preserves a real analytic volume form. Assume that $X$ has no periodic orbits, then $M$ is either diffeomorphic to the three dimensional torus $\ttt^3$ or to a parabolic torus bundle over the circle.
\end{Teorema}

\begin{Teorema}\label{teoremaperiogeodana}
Assume that $X$ is a geodesible vector field on an oriented closed 3-manifold $M$ that is not a torus bundle over the circle. Then if $X$ is real analytic, it possesses a periodic orbit.
\end{Teorema}

A parabolic torus bundle over the circle is a torus bundle over the circle where the attaching map is isotopic to a parabolic element of $SL(2,\zz)$. We will like to point out that every parabolic torus bundle over the circle admits a aperiodic $C^\omega$ geodesible vector field: consider the suspension of the parabolic matrix composed with an irrational translation.

Before explaining the main tools we will use in this article, let us briefly review the results concerning the existence of periodic orbits for non singular vector fields on closed 3-manifolds. Let $X$ be a $C^\infty$ non singular vector field on an oriented closed 3-manifold $M$, we have the next results

\vspace{3mm}

{\small
\tracingtabularx
\begin{tabularx}{400pt}{|c|X|X|}
\hline
$X$ & &  {\bf Volume preserving} \\ \hline
 & (1) On any $M$ there exist examples of $C^\infty$, or even $C^\omega$, vector fields without periodic orbits. & (2) On any $M$ there exist examples of $C^1$ vector fields without periodic orbits.  \\ \hline
{\bf Geodesible} & (3) Theorem \ref{teoremaperiogeodana}. & (4) If $M$ is not a torus bundle over the circle, $X$ has a periodic orbit. \hfill ($4^\prime$) If $X$ is $C^\omega$, we have theorem \ref{teoremaperiogeodanavol}.
\\ \hline
{\bf Reeb} & & (5) On any $M$, the vector field has periodic orbits.\\ \hline
\end{tabularx}
}

\vspace{3mm}

Theorems (1) and (2) where proved using plugs. Plugs where first used by W. Wilson \cite{wils} and are used to change a vector field inside a flow-box in order to break periodic orbits. The construction of the examples in (1) is due to K. Kuperberg \cite{kkup}, and in (2) to G. Kuperberg \cite{gkup}. We will like to point out that D. Sullivan's classification (\cite{sull3}) of geodesible vector fields implies that plugs cannot be used for geodesible vector fields (we refer the reader to \cite{these}). Result number (5) was proved for $M$ diffeomorphic to the three sphere $\sss^3$ or for $M$ with non trivial second homotopy group by H. Hofer \cite{hofe2}, and for any 3-manifold by C. H. Taubes \cite{taub}. These theorems are valid for $C^2$ vector fields.

As we said, Reeb vector fields and vector fields that admit a global cross section are geodesible. Hence, for the first type of vector fields theorems \ref{teoremaperiogeodanavol} and \ref{teoremaperiogeodana} are a consequence of (5), and for the second type of vector fields we just need to use Lefschetz fixed point formula. As we will explain, there are geodesible vector fields that are neither Reeb or suspensions. For these we will find an invariant set, in whose complement the vector field is a Reeb vector field. The real analyticity of this invariant set will allow us to study its complement and find a periodic orbit in it.

 Finally we will like to point out some facts about (4). During the preparation of this work, M. Hutchings and C. H. Taubes announced the result corresponding to (4) with a slightly different definition: if $M$ is an oriented closed 3-manifold with a stable Hamiltonian structure, and $X$ is  the associated Reeb vector field, then either $X$ has a periodic orbit or $M$ is a torus bundle over the circle \cite{huta}. This is equivalent to (4). A stable Hamiltonian structure on $M$ is a pair $(\alpha, \omega)$, where $\alpha$ is a 1-form and $\omega$ is a closed 2-form on $M$, such that
\begin{eqnarray*}
\alpha\wedge\omega & > & 0\\
\ker(\omega) & \subset & \ker(d\alpha).
\end{eqnarray*}
Note that the first condition implies that $\omega$ is non vanishing, and consequently $d\alpha=f\omega$ for a function $f$ on $M$. The Reeb vector field $X$ of a stable Hamiltonian structure is determined by the equations
$$
\alpha(X)=1 \qquad \mbox{and} \qquad \iota_X\omega=0.
$$
Thus, a volume preserving geodesible vector field is the Reeb vector field of a stable Hamiltonian structure, just put $\omega=\iota_X\mu$ where $\mu$ is the invariant volume form. Conversely, a Reeb vector field of a stable Hamiltonian structure is a geodesible volume preserving vector field, the volume is given by $\alpha\wedge \omega$.

In their proof they use the following theorem, that is proved using Seiberg-Witten invariants (we refer to theorem {\bf 1.2} of \cite{huta}). We say that a periodic orbit is nondegenerate if the eigenvalues of Poincar\'e's first return map are different from 1, and that it is elliptic if the eigenvalues are in the unit circle.

\begin{teorema}[Hutchings, Taubes]
Let $M$ be a closed oriented 3-manifold with a contact form $\alpha$ such that all the periodic orbits of its Reeb vector field are nondegenerate and elliptic. Then $M$ is a lens space, there are exactly two periodic orbits and they are the core circles in the solid tori of a genus one Haegaard splitting of $M$.
\end{teorema}

Let $T$ denote a 2-torus, a consequence of this theorem is (we refer to section {\bf 5.2} of \cite{huta}) 

\begin{teorema}[Hutchings, Taubes]
\label{teohuta}
Let $N$ be a compact 3-manifold with boundary that is endowed with a Reeb vector field $X$ (of a contact form) tangent to $\partial N$. Then if $N$ is not diffeomorphic to $T\times [0,1]$, the vector field $X$ possesses a periodic orbit.
\end{teorema}

We will use this theorem when proving theorem \ref{teoremaperiogeodana}. As we will explain in section \ref{extra}, we can prove the following two theorems without using their theorem, and without using Seiberg-Witten invariants. The second one is a particular case of the main result in \cite{huta}.

\begin{Teorema}\label{teoremaperiogeodanadebil}
Assume that $X$ is a geodesible vector field on an oriented closed 3-manifold $M$, that is either diffeomorphic to $\sss^3$ or has non trivial $\pi_2$. Then if $X$ is real analytic, it possesses a periodic orbit.
\end{Teorema}

\begin{Teorema}\label{teoremaperiogeodvol}
Assume that $X$ is a geodesible vector field on an oriented closed 3-manifold $M$, that is either diffeomorphic to $\sss^3$ or has non trivial $\pi_2$. Then if $X$ is of class $C^2$  and preserves a volume, it possesses a periodic orbit.
\end{Teorema}

The text is organized as follows, we will begin by giving examples of geodesible vector fields. In section \ref{rotacional} we will give another characterization of these vector fields on 3-manifolds in terms of their curl. The proof of  theorem \ref{teoremaperiogeodanavol} is given in section \ref{preservavolumen},  and the one of theorem \ref{teoremaperiogeodana} in section \ref{periodicageod}.

I will like to thank my Ph.D. thesis advisor \'Etienne Ghys for all his patience and comments during the preparation of this work. I will also thank Mexico, that has sponsored my Ph.D. studies through the scholarship program of the Consejo Nacional de Ciencia y Tecnolog\'ia (CONACyT). Special thanks to Patrice Le Calvez, Pierre Py and Klaus Niederkr\"uger for all the comments and  long discussions. 

\section{Geodesible vector fields}
\label{geodesibles} 

Geodesible vector fields were first studied by H. Gluck. He was interested in filling manifolds with geodesics, in other words which manifolds can be foliated by geodesics \cite{gluc}.

Let $X$ be a geodesible vector field and $g$ a Riemannian metric making the orbits of $X$ geodesics. We can consider the differential 1-form $\alpha=\iota_Xg$. We can reparameterize $X$ so that $\alpha(X)=1$. The condition $\nabla_XX=0$ where $\nabla$ is the Levi-Civita connection, {\it i.e.} the orbits are geodesics for $g$, is equivalent to $L_X\alpha=0$. Then $X\in\ker(d\alpha)$. The plane field defined by $\xi=\ker(\alpha)$ is then invariant.
Conversely, the existence of a 1-form such that
$$
\alpha(X)=1\qquad \mbox{and} \qquad X\in\ker(d\alpha),
$$
implies that $X$ is geodesible, we just need to take a Riemannian metric making $X$ of unit length and orthogonal to $\ker(\alpha)$. 

Assume that the ambient manifold $M$ has dimension three. If the plane field $\xi=\{\ker(\alpha)\}$ is a contact structure, or equivalently $\alpha\wedge d\alpha\neq 0$, $X$ is the Reeb vector field of $\alpha$. On the other hand, if $\alpha\wedge d\alpha=0$ the form $\alpha$ is closed, implying by S. Schwartzman's \cite{schwa} and D. Tischler's theorems  \cite{tisc} that $M$ is a fiber bundle over $\sss^1$ and each fiber is a global section of $X$: a closed submanifold that intersects every orbit. Thus geodesible vector fields contain Reeb vector fields associated to a contact structure and vector fields that admit a global cross section.

Let us mention some examples of geodesible flows. 
\begin{list}
{\arabic{num}.}{\usecounter{num}}
\item {\it Geodesic flow on a Riemannian surface.} Let $S$ be a Riemannian surface. The vector field generating the geodesic flow is in the unit tangent bundle $T^1S$. The orbits in this space are geodesics for {\it Sasaki metric}. Moreover, it is the Reeb vector field of the 1-form given by the contraction of this metric with the vector field. Using Stoke's theorem we can prove that Reeb vector fields do not admit global cross sections, and thus the vector field generating the geodesible flow on $S$ does not admit cross sections. For the details we refer to G. Paternain's book \cite{pate}.

\item {\it Killing vector fields.} Let $M$ be a Riemannian manifold, denote the metric by $g$. A vector field $X$ on $M$ which generates an isometric flow on $M$ is a Killing vector field. Such a vector field is geodesible: there exists a Riemannian metric $g^\prime$ on $M$, conformal to $g$, for which the trajectories of $X$ are geodesics. 

\item{\it Vector fields tangent to 1-foliations by closed curves.} Consider a manifold $M$ foliated by closed curves. A. W. Wadsley showed that there is a Riemannian metric on $M$ making these curves geodesics if and only if their lengths are bounded, \cite{wads}. D. Epstein had previously showed that this bounded length condition is always satisfied in dimension three, \cite{epst}. In \cite{sull2}, Sullivan showed that it can fail in dimension strictly greater than four, and Epstein and E. Vogt constructed a flow in a compact 4-manifold such that all the orbits are circles with unbounded length \cite{epvo}. Thus in dimension three a one dimensional foliation by closed curves can be seen as the orbits of a geodesible vector field.

In dimension three, these examples are contained in the precedent ones. There is a vector field tangent to the foliation by circles whose flow is given by a locally free action of $\sss^1$, see \cite{epst}. Such a flow is isometric, the invariant metric is obtained as follows: beginning with any metric on $M$, we take the mean of the transformations of the metric under the $\sss^1$ action, relatively to the Haar measure on $\sss^1$.
\end{list}

\section{Geodesible vector fields are Beltrami}
\label{rotacional}

In this section we will consider geodesible vector fields on closed 3-manifolds. We will  establish a correspondence between geodesible vector fields and some solutions of the Euler equation of an ideal stationary fluid, known in hydrodynamics as Beltrami vector fields. 
Our result is a generalization of theorem {\bf 2.1} of \cite{etghi} by J. Etnyre and R. Ghrist: {\it the class of Reeb vector fields of a contact structure on a 3-manifold is identical to the class of vector fields that  have non zero curl and are colinear with it, for a suitable Riemannian metric.}

\begin{proposicion}\hspace{2mm}
\label{teorelgeodcurl}
Let $M$ be an oriented 3-manifold. Any $C^2$ vector field that is parallel to its curl, for a Riemannian metric, is geodesible. Conversely, any geodesible vector field, modulo a reparameterization, is parallel to its curl.
\end{proposicion}

The definition of the curl of a vector field in $\rr^3$ depends upon a Riemannian metric. We adopt the following definition: the curl of a vector field $X$ on a Riemannian 3-manifold $M$, with metric $g$ and arbitrary distinguished volume form $\mu$, is the unique vector field $\mbox{curl}(X)$ given by
$$
\iota_{\mbox{curl}(X)}\mu=d\iota_Xg.
$$
Taking the curl with respect to an arbitrary volume form makes the subsequent results valid for a more general class of fluids: for example {\it basotropic flows}, these are compressible for the Riemannian volume and incompressible for a rescaled volume form. We refer the reader to section {\bf VI.2.A} of V. I. Arnold and B. Khesin's book \cite{arkh}.
When $\mu$ is the Riemannian volume the curl assumes the more common form
$$
\mbox{curl}(X)=\psi(*d\iota_Xg),
$$
where $*$ is the Hodge star operator, and $\psi$ is the isomorphism between vector fields and differential 1-forms derived from $g$.

\begin{definicion}
The Euler equation of an ideal incompressible fluid on a Riemannian manifold $M$ endowed with a volume form $\mu$, is 
\begin{eqnarray*}
\frac{\partial X_t}{\partial t} & = & -\nabla_{X_t}X_t-\mbox{grad}(p)\\
L_{X_t}\mu & = & 0,
\end{eqnarray*}
where the velocity vector field $X_t$ and the function $p$ are time dependent. The second equation means that $X_t$ preserves the volume form $\mu$. 
\end{definicion}

We will deal with  the Euler equation of an ideal steady fluid on $M$. That is, the vector field $X$ will be time independent and not necessarily volume preserving. We get the equation $\nabla_XX=-\mbox{grad}(p)$,
for a pressure function $p$. Using the identity $\nabla_XX=X\times\mbox{curl}(X)+\frac{1}{2}\|X\|$, from page {\bf 588} of \cite{abmr}, we can reduce the equation to the form
\begin{equation}\label{euler}
X\times \mbox{curl}(X)=\mbox{grad}(b)
\end{equation}
where $b=-p-\frac{1}{2}\|X\|$. The function $b$ is known as the {\it Bernoulli function} of $X$.

\vspace{1mm}

\sdemostracion{\textbf{of proposition \ref{teorelgeodcurl}. \hspace{2mm}} Assume first that $X$ is a geodesible vector field. We know that, modulo a reparameterization, $\iota_Xg=\alpha$ is an invariant 1-form and that $X\in \ker(d\alpha)$. Using the definition of the curl we have that $\iota_{\mbox{curl}(X)}\mu=d\alpha$, thus $X\in \ker(\iota_{\mbox{curl}(X)}\mu)$. Since $\mu$ is a volume form we have that $\mbox{curl}(X)=fX$, for a function $f:M\to\rr$. Observe that this function can be zero. 

Conversely, if a vector field $X$ is such that $\mbox{curl}(X)=fX$ with respect to a Riemannian metric $g$, setting $\iota_Xg=\alpha$ we have that $\alpha(X)>0$ and $\iota_Xd\alpha=\iota_X\iota_{fX}\mu=0$. We can rescale $X$ so that $\alpha(X)=1$. Thus the vector field $X$ is geodesible.
}

\begin{definicion}
A vector field such that $\mbox{curl}(X)=fX$ for a function $f$ on $M$ is a Beltrami vector field in hydrodynamics. In magnetodynamics these vector fields are known as {\it force-free} vector fields. 
\end{definicion}

Before finishing this section let us analyze volume preserving geodesible vector fields. Let $X$ be a geodesible vector field preserving a volume form $\mu$. 
An important consequence of the results above is that the function $f$ is constant along the orbits of $X$. This follows from
\begin{equation}\label{rotar}
0=L_{\mbox{curl}(X)}\mu=L_{fX}\mu=fd\iota_X\mu+df\wedge\iota_X\mu.
\end{equation}
Since $d\iota_X\mu=0$, we have that $f$ is a first integral of $X$. 

When $f$ is different from zero, $X$ is a Reeb vector field of the contact form $\alpha$. We can choose the volume form $\mu=\alpha\wedge d\alpha$ and thus $f=1$. Examples coming from hydrodynamics are the $ABC$ vector fields on the three torus $\ttt^3$. On the three dimensional torus, $\{(x,y,z)| \mbox{mod}\,2\pi\}$, an $ABC$ vector field is defined by the equations
\begin{eqnarray*}
v_x & = & A\sin z+C\cos y\\
v_y & = & B\sin x+A\cos z\\
v_z & = & C\sin y+B\cos x.
\end{eqnarray*}
They preserve the unit volume form and have $\mbox{curl}(v)=v$. These vector fields where first studied by I. S. Gromeka in 1881, rediscovered by E. Beltrami in 1889, and largely studied in the context of hydrodynamics during the last century. When one of the parameters $A$, $B$ or $C$ vanishes, the vector field is integrable.  By symmetry of the parameters, we may assume that $1=A\geq B\geq C\geq 0$. 

In 1986, T. Dombre, U. Frisch, J. Greene, M. H\'enon, A. Mehr and A.  Soward \cite{dfgh}, showed the absence of integrability when $ABC\neq 0$. They also showed that under the precedent convention, the vector field is non singular if and only if $B^2+C^2<1$. Though the list of publications concerning $ABC$ vector fields is extensive, there is very little known about the global features of these vector fields, apart from cases when $C$ is zero or a perturbation thereof.

\section{Aperiodic volume preserving $C^\omega$ vector fields}
\label{preservavolumen}

In this section we will prove theorem \ref{teoremaperiogeodanavol}, thus we will classify the 3-manifolds that admit aperiodic volume preserving $C^\omega$ vector fields on closed 3-manifolds.

\sdemostracion{{\bf of theorem \ref{teoremaperiogeodanavol}}{\hspace{2mm}} We know, from section \ref{rotacional} that $\mbox{curl}(X)=fX$ for a real analytic function $f:M\to\rr$. Further, since $X$ preserves a volume given by a real analytic differential form that we will call $\mu$, the function $f$ is a first integral of $X$ (see equation \ref{rotar}). Let $\xi$ be the plane field defined by the kernel of $\alpha$. We can distinguish the next three cases:
\begin{itemize}
\item[(I)] $f$ is never zero. In this case the plane field $\xi$ is a contact structure and $X$ is the associated Reeb vector field. Thus H. Hofer's and C. H. Taubes' theorems imply the existence of a periodic orbit of the vector field $X$ on any oriented closed \mbox{3-manifold}. Thus under the hypothesis of the theorem, $X$ cannot be a Reeb vector field.

\item[(II)] $f$ is identically zero. As we previously said, this implies that the 1-form $\alpha$ is closed. Thus the vector field has a section $T$, an oriented closed surface without boundary. Since $X$ is aperiodic, the Lefschetz fixed point theorem (we refer to theorem {\bf 8.6.2} of \cite{kaha}) implies that $T$ must be a torus and $M$ is a torus bundle over the circle.

Consider an aperiodic geodesible vector field, then
$$
M ={T\times [0,1]}/{(x,0)\sim (h(x),1)},
$$
where $h$ is the first return map defined by $X$. Let $T_0=T\times \{0\}$. The diffeomorphism type is defined by the isotopic class of $h$. We can define a function $g:T_0\to \rr$ that is the first return time: for a point $x\in T_0$, the value of $g(x)$ is the time the orbit of $x$ takes to return to $T_0$. Reparameterizing the flow of $X$ with $g$, we get that it is $C^\omega$ conjugated to the suspension of the diffeomorphism $h$. 

In order to characterize the diffeomorphism type of $M$, we will use the Lefschetz fixed point formula to prove.

\begin{lema}
The trace of $h$ in $H^1(T)\simeq \zz^2$ is equal to two.
\end{lema}

\demostracion{By Lefschetz fixed point formula the alternated sum of the trace of $h$ in $H^i(T)$, for $i=0,1,2$ is zero. Clearly, the trace in $H^0(T)$ and in $H^2(T)$ is equal to one, thus proving the lemma. 
}

We conclude that $M$ is either diffeomorphic to $\ttt^3$ or a parabolic torus bundle over the circle.

\item[(III)] $f$ is equal to zero on a compact invariant set $f^{-1}(0)=A\subset M$. As we previously said $A$ is the set where $\alpha$ is closed. Observe that for a regular value $a$ of $f$, the compact set $f^{-1}(a)$ is a finite union of disjoint invariant tori. 

Let us study the topology of the critical levels. Consider now a critical value $c$ of $f$.

\begin{lema}\label{Aanalitico}
Each connected component $C$ of $f^{-1}(c)$ is homeomorphic to a torus and $X|_C$ is topologically conjugate to a linear irrational vector field.
\end{lema}

\demostracion{Since $f$ is a real analytic function, $C$ is a real analytic set, thus it is a Whitney stratified set: it can be decomposed into manifolds of dimension less or equal to two. Take $x\in C$.  The non singularity of the vector field $X$ yields to the existence of a flow box $N\simeq D\times [-1,1]$ with $D$ a transverse disc. Assume that $x\in D_0=D \times \{0\}$. Since $C$ is invariant under the flow induced by $X$ we have that
$$
N\cap C\simeq (D\cap C)\times [-1,1].
$$
We know that the dimension of the strata manifolds that compose $C$ is at most two. Using H. Whitney's theory, we get that $D_0\cap C$ is homeomorphic to a radial $k$-tree centered at $x$. This $k$-tree is invariant under the flow. We refer the reader to H. Whitney's book \cite{whit}.

If $k=0$, the set $C$ is of dimension one and compact, thus it is a periodic orbit, a contradiction to the hypothesis of the lemma. If $k>2$, the point $x$ is contained in a dimension one submanifold of $C$. This submanifold is a periodic orbit, a contradiction. Then for any $x\in C$, the intersection $D_0(x)\cap C$ is homeomorphic to an invariant 2-tree, thus $C$ is an invariant surface. The argument we used in the case $k>2$, implies also that the 2-tree is a non singular $C^\omega$ curve in the disc, and thus $C$ is a non singular real analytic oriented surface that admits a non singular vector field. Then $C$ is a torus. Since $X$ has no periodic orbits on $C$ it is topologically conjugate to a linear irrational vector field.
}

Then all the levels of $f$ are invariant tori, thus we have decomposed the manifold into tori.

{\bf Explicit expression for $X$ in $M$}

On each torus there is a tangent non singular vector field $Y$ defined by the equation $\alpha(Y)=0$ and $X\times Y=N$ where $N$ is the normal vector field to the torus.  The reason why it is non singular is that the 1-foliation of each torus that is tangent to $\xi$ is non singular. The vector fields $X$ and $Y$ are linearly independent and commute. 

Take one of the tori and call it $T_0$, then there is a neighborhood $U$ of it that is diffeomorphic to $T\times[0,1]$ where $X$ is tangent to $T_t$, with $T_t=T\times\{t\}$. We claim that we can give an explicit expression for $X$ in $U$.

\begin{lema}
There are $C^\infty$ functions $a_1, a_2, a_3, a_4$ defined on $[0,1]$ such that the vector fields
$$
a_1(t)X+a_2(t)Y \qquad \mbox{and} \qquad a_3(t)X+a_4(t)Y
$$
are linearly independent on $T_t$ and all their orbits are periodic of period one.
\end{lema}

\demostracion{Fix $t\in[0,\delta]$. Denote by $\phi_s$ the flow of $X$ and $\psi_s$ the flow of $Y$ on $T_t$. For a fixed point $x\in T_t$, consider the map
\begin{eqnarray*}
\Phi :\rr^2 & \to & T_t\\
(s_1,s_2) & \mapsto & \phi_{s_1}\psi_{s_2}(x).
\end{eqnarray*}
Since $X$ and $Y$ are linearly independent and commute, $\Phi$ is a covering map. Then for $y\in T_t$, the inverse image $\Phi^{-1}(y)$ is a lattice in $\rr^2$. 

Define the functions $a_i$ for $i=1,2,3,4$ such that $(a_1(t), a_2(t))$ and $(a_3(t), a_4(t))$ form a basis for the lattice for each $t$.  Then the vector fields
$$
a_1(t)X+a_2(t)Y \qquad \mbox{and} \qquad a_3(t)X+a_4(t)Y
$$
are linearly independent and have closed orbits of period one.
}

Hence in $U$ we have a system of coordinates $(x,y,t)$ such that $f(x,y,t)=f(x,y,0)+t$, and $X$ is a linear flow on each torus $T_t$ that can be written as 
$$
X=\tau_1(t)\frac{\partial}{\partial x}+\tau_2(t)\frac{\partial}{\partial y}.
$$
Since $X$ is aperiodic the ratio $\frac{\tau_1(t)}{\tau_2(t)}$ is constant and equal to an irrational number. We claim that $\tau_1$ and $\tau_2$ are independent of $t$. Since $L_X\alpha=0$ we have that 
$$
L_X\left(\alpha\left(\frac{\partial}{\partial t}\right)\right)=-\alpha\left(\tau_1^\prime(t)\frac{\partial}{\partial x}+\tau_2^\prime(t)\frac{\partial}{\partial y}\right) =-\frac{\tau_1^\prime(t)}{\tau_1(t)}.
$$
The right side of the equation depends only on $t$ and is constant on each torus $T_t$. Moreover, it is a coboundary, so it is zero on each torus. This implies that $\tau_1$ is constant, and since $\frac{\tau_1}{\tau_2(t)}$ is constant, the function $\tau_2$ is constant. Thus,
\begin{equation}\label{X}
X=\tau_1\frac{\partial}{\partial x}+\tau_2\frac{\partial}{\partial y}.
\end{equation}

This implies that locally $f$ behaves as a projection into an interval. Then $M$ is a torus bundle over a compact 1-manifold, {\it i.e} it is a torus bundle over the circle and $X$ is tangent to each fiber. 

To finish the proof of the theorem we have to classify these manifolds. Let us cut $M$ along a torus, we get $T\times [0,1]$ endowed with the vector field given in equation \ref{X}.

{\bf Explicit expression for $\alpha$}

The 1-form $\alpha$ can be written as
$$
\alpha=A_1dx+A_2dy+A_3dt,
$$
where $A_1$, $A_2$ and $A_3$ are functions of $(x,y,t)$. Using the fact that $\alpha(X)=1$ we get that $A_2=\frac{1-\tau_1A_1}{\tau_2}$. Then,
\begin{equation*}
d\alpha=\left(\frac{\partial A_3}{\partial x}-\frac{\partial A_1}{\partial t}\right)dx\wedge dt+\left(\frac{\partial A_3}{\partial y}+\frac{\tau_1}{\tau_2}\frac{\partial A_1}{\partial t}\right)dy\wedge dt-\left(\frac{\tau_1}{\tau_2}\frac{\partial A_1}{\partial x}+\frac{\partial A_1}{\partial y}\right)dx\wedge dy.
\end{equation*}
Looking at the third term in the expression of $\iota_Xd\alpha$ and using that $X\in\ker(d\alpha)$, we get that $A_1$ and $A_2$ are functions of $t$. Moreover, we can put
$$
A_1(t)=\gamma(t)\tau_2+c \qquad \mbox{and} \qquad A_2(t)=-\gamma(t)\tau_1+\frac{1-c\tau_1}{\tau_2},
$$
where we added the second term in both expressions to satisfy the condition $\alpha(X)=1$, with $c$ any non zero constant. Looking at the first two terms in the expression of $\iota_Xd\alpha$, we obtain
that $A_3$ is a function of $t$ and hence
$$
\alpha=\gamma(t)(\tau_2dx-\tau_1dy)+A_3(t)dt+\frac{c\tau_2dx+(1-c\tau_1)dy}{\tau_2}.
$$

We have that $M$ is $T\times[0,1]$ whose boundary components are identified with a diffeomorphism $h$ of the torus. We know that $h$ preserves the linear foliation given by the kernel of $\alpha$ restricted to the boundary tori. The slope depends on $\gamma(t)$ and the constant $c$. We can choose $T_0$ in such a way that this 1-foliation has irrational slope. Moreover, $h$ preserves the vector field $X$. Let $\Phi$ be the induced map on the universal cover of the torus $\rr^2$, we have that
$$
\Phi(x,y)=(g_1(x,y),g_2(x,y)).
$$

\begin{lema}
The functions $g_1$ and $g_2$ are linear in $x$ and $y$.
\end{lema}

\demostracion{Let us call $\delta$ the slope of the foliation of the boundary tori in $T\times[0,1]$ defined by the kernel of $\alpha$. As we said, we took $T$ such that $\delta$ is irrational. Then there is a function $G$ such that
$$
g_2(x,y)-\delta g_1(x,y)=G(y-\delta x).
$$
Moreover, the functions $g_1$ and $g_2$ can be expressed in a unique form as the sum of a linear function and a periodic function:
$$
g_1=l_1+p_1 \qquad \mbox{and}\qquad g_2=l_2+p_2.
$$
Then, $(l_2-\delta l_1)(x,\delta x)+(p_2-\delta p_1)(x,\delta x)=G(0)$. The function $p_2-\delta p_1$ is periodic, hence bounded. The function
$$
x\mapsto  (l_2-\delta l_1)(x,\delta x),
$$
is linear and bounded, hence it is equal to zero. This implies that $p_2-\delta p_1$ is constant on a straight line with irrational slope, then it is constant everywhere. Hence,
$$
G(y)=g_2(0,y)-\delta g_1(0,y)=(l_2-\delta l_1)(0,y) +(p_2-\delta p_1)(0,y)=ay +b,
$$
for two real numbers $a,b$. Hence $g_2(x,y)-\delta g_1(x,y)=a(y-\delta x)+b$. 

Since $h$ preserves $X$, it preserves the 1-foliation given by its orbits. Thus there exist real numbers $c,d$ such that
$$
g_2(x,y)-\frac{\tau_2}{\tau_1}g_1(x,y)=c\left(y-\frac{\tau_2}{\tau_1}x\right)+d.
$$
We get that $\delta g_1(x,y)+a(y-\delta x)+b= \frac{\tau_2}{\tau_1}g_1(x,y)+c(y-\frac{\tau_2}{\tau_1}x)+d$. We conclude that $g_1$ and $g_2$ are linear functions on $x$ and $y$.
}

Then $h\in SL(2,\zz)$ and one of its eigenvalues must be one. Thus $h$ is either conjugate to the identity or to a parabolic matrix.

\end{itemize}
}

\section{Existence of periodic orbits for $C^\omega$ geodesible vector fields}
\label{periodicageod}

\sdemostracion{{\bf of theorem \ref{teoremaperiogeodana}}{\hspace{2mm}} Take a real analytic volume form on the ambient manifold, and let $f$ be the real analytic function satisfying $\mbox{curl}(X)=fX$. The set $f^1(0)$ does not depends on the choice of the volume form. Observe that for the cases (I) and (II) in the precedent proof we did not use the volume preserving hypothesis, hence we only need to prove the theorem in case (III).  

Denote by $\xi$ the plane field defined by the kernel of the differential 1-form $\alpha=\iota_Xg$, where $g$ is the Riemannian metric making the orbits of $X$ geodesics.
\begin{itemize}
\item[(III)] The function $f$ is equal to zero on a compact invariant set $A\subset M$. Moreover, $A$ is real analytic. Assume that $X$ does not have any periodic orbit on $A$, by lemma \ref{Aanalitico} the set $A$ is a finite union of invariant tori. Let $B$ be a connected component of $M\setminus A$, and assume that $\alpha\wedge d\alpha\geq 0$ in $\overline{B}$. The case where $\alpha\wedge d\alpha\leq 0$ being equivalent. Let us consider $\overline{B}\subset M$.

The idea of the rest of the proof is to approximate the plane field $\xi$ on $\overline{B}$ by a contact structure, that is transverse to $X$. The vector field $X$ will be a Reeb vector field of the new contact structure. Then by theorem \ref{teohuta} we conclude that $X$ possesses a periodic orbit in $\overline{B}$. We will use the following proposition. 

\begin{proposicion}\label{propaproximacion}
The $C^\infty$ plane field $\xi$ on $\overline{B}$  can be $C^\infty$ approximated  by a contact structure $\eta$. Moreover, a Reeb vector field of $\eta$ is $X$.
\end{proposicion}

To prove the proposition we will begin by a lemma. Consider $\rr^3$ with coordinates $(x,y,t)$.

\begin{lema}\label{lemaaproximacion}
Let $\xi$ be a $C^k$ plane field on
$$
V=\{|x|\leq 1, |y|\leq 1, 0\leq t\leq 1\}\subset \rr^3,
$$
given by the kernel of the 1-form $\beta=dx-a(x,y,t)dy$ such that $\beta\wedge d\beta\geq 0$. Suppose that the plane field is a contact structure near \mbox{$\{t=1\}$.} Then $\xi$ can be approximated by a plane field $\xi^\prime$ which coincides with $\xi$ together with all its derivatives along $\partial V$ and is contact inside $V$.
\end{lema}

\demostracion{Observe that $\xi$ is transversal to the $x$-curves and tangent to the $t$-curves. We have that
$$
\beta\wedge d\beta=\frac{\partial a}{\partial t}(x,y,t)dx\wedge dy\wedge dt\geq 0,
$$
then $\frac{\partial a}{\partial t}(x,y,t)\geq 0$ in $V$ and $\frac{\partial a}{\partial t}(x,y,1)>0$. Then there exists a function $\tilde{a}(x,y,t)$ such that
$$
\frac{\partial \tilde{a}}{\partial t}(x,y,t)>0
$$
in the interior of $V$ and coincides with $a$ along $\partial V$. Moreover, along $\partial V$ all the derivatives of $a$ and $\tilde{a}$ coincide. Then the plane field
$$
\xi^\prime=\{\ker(dx-\tilde{a}(x,y,t)dy)\}
$$
is the perturbation with the required properties.
}

\sdemostracion{{\bf of the proposition.}\hspace{2mm} Let $\xi$ be the positive plane field of $\overline{B}$ of class $C^\infty$ and $C(\xi)=B$ its contact part. For every point $p\in \partial \overline{B}$, choose a simple curve $\gamma_p$ which is tangent to $\xi$, begins at $p$ and ends at a point $p^\prime\in C(\xi)$. Let $S_p:[0,1]\times [0,1]\to M$ be an embedding such that the image of $[0,1]\times\{\frac{1}{2}\}$ is $\gamma_p$ and gives us a surface in $M$ that is transverse to $X$, and such that $S_p([0,1]\times\{\cdot\})$ is tangent to $\xi$. Moreover, we will ask that the image of $\{0\}\times[0,1]$ is contained in $\partial\overline{B}$.

Observe that since the orbits of $X$ are geodesics its flow $\phi_s$ preserves orthogonality. Then the images under the flow of the curves $S_p([0,1]\times\{\cdot\})$ are tangent to $\xi$. Pushing the above surface with the flow we get a region
$$
\mathcal{V}_p=\{\phi_s\cdot S_p([0,1]\times[0,1])\,|\, s\in[-\epsilon, \epsilon]\},
$$
for a given $\epsilon$. Observe that part of the boundary of this region is in $\partial\overline{B}$ and the surface $\phi_s\cdot S_p([0,1]\times[0,1])$ is transverse to $X$ for every $s$.

Denote $V=[-1,1]\times[-1,1]\times[0,1]$ with coordinates $(x,y,t)$. Then, there exists an embedding
$$
F_p:V\to \mathcal{V}_p\subset\overline{B}
$$
satisfying
\begin{itemize}
\item[(i)] the line segment $(0,0,t)$, with $t\in[0,1]$, is mapped to $\gamma_p$;
\item[(ii)] the images of the $t$-curves are tangent to $\xi$;
\item[(iii)] the the vector field $\frac{\partial}{\partial x}$ is mapped to $X$;
\item[(iv)] the image of $(x,y,1)$ is in $C(\xi)$ for all pairs $(x, y)$;
\item[(v)] the image of $(x,y,0)$ is in $\partial\overline{B}$ for all pairs $(x, y)$.
\end{itemize}

Let $W$ be the interior of $[-\frac{1}{2},\frac{1}{2}]\times[-\frac{1}{2},\frac{1}{2}]\times[0,1]$ and $W^\prime$ the interior of $[-\frac{3}{4},\frac{3}{4}]\times[-\frac{3}{4},\frac{3}{4}]\times[0,1]$. In the manifold $\overline{B}$ we can find a finite number of points $p_1,p_2,\ldots ,p_n$ and corresponding paths, such that the open sets $W_i=F_{p_i}(W)$ cover an open neighborhood of $\partial \overline{B}$. Set 
$$
W_i^\prime=F_{p_i}(W^\prime) \qquad \mbox{and} \qquad V_i=F_{p_i}(V).
$$
Thus we have that $W_i\subset W_i^\prime\subset V_i$ for every $i$. The pull back $\beta_i=(F_{p_i})^*\alpha$ is a 1-form such that $\frac{\partial}{\partial t}$ is in its kernel, $\beta_i(\frac{\partial}{\partial x})=1$ and $\frac{\partial}{\partial x}$ is in the kernel of $d\beta_i$. Thus we can write 
$$
\beta_i=dx-a_i(x,y,t)dy,
$$
The condition $\frac{\partial}{\partial x}\in \ker(d\beta_i)$ implies that the function $a_i$ must be independent of $x$. Hence, $\beta_i=dx-a_i(y,t)dy$ and $\frac{\partial a_i}{\partial y}=0$ only on the surface $\{t=0\}$.
 
Applying the lemma, we can perturb $\xi_1$ into a plane field $\xi_1^\prime$ which is contact in $W^\prime$. Since $\xi_1^\prime$ is defined by the kernel of the differential form $dx-\tilde{a}_1(y,t)dy$. The latter form coincides with $\beta_i^\prime$ along the  $\partial V\setminus \{t=0\}$. The vector field $X$ is in the kernel of $\frac{\partial \tilde{a}_1}{\partial t}(y,t)dy\wedge dt$.  The push forward of $\xi_1^\prime$ defines a perturbation of $\xi$ in $\overline{B}$ that is contact in $W_1^\prime$, and in $W_1^\prime$ a Reeb vector field of the new plane field is $X$.

Unfortunately, we cannot simply continue the process because the perturbation inside $V_1$ affects the properties of the rest of the embeddings. However, if the perturbation is small enough, it is possible to modify the embeddings, for $i=2,3,\ldots,n$, into $F_{p_i}^\prime$ satisfying conditions (ii) through (iv), and the condition
$$
W_i\subset F_{p_i}^\prime(W^\prime) \qquad \mbox{for} \qquad i=2,3,\ldots, n,
$$
at the place of (i). This is sufficient to continue with the process of perturbation. Since the $W_i$ cover a neighborhood of  $\partial \overline{B}$, we get a contact structure $\eta$ in $\overline{B}$ such that one of its Reeb vector fields is parallel to $X$.
}

Hence we have a contact structure $\eta$ on $\overline{B}$ that has a Reeb vector field that is parallel to $X$. Then by theorem \ref{teohuta} we have a periodic orbit of $X$ in $\overline{B}$. This finishes the proof.
\end{itemize}
}

 \section{Proofs without Seiberg-Witten invariants}
\label{extra}

In this section we will like to give a proof of theorems \ref{teoremaperiogeodanadebil} and \ref{teoremaperiogeodvol} without using theorem \ref{teohuta}. We will begin by explaining the main result we use and its proof. The theorem was proved in the case of solid tori by R. Ghrist and J.Etnyre \cite{etghII}. 

\begin{teorema}\label{teoremafrontera}	
Let $N$ be a 3-manifold with boundary endowed with a Reeb vector field $X$ tangent to the boundary. If $N$ is either diffeomorphic to a solid torus or $\pi_2(N)\neq 0$, the vector field $X$ possesses a periodic orbit.
\end{teorema}

Let us start by introducing some concepts of contact geometry. The vector field $X$ is the Reeb vector field of a contact form $\alpha$ defined on $N$, and we will call $\xi$ the contact structure. Given an embedded surface $S$ in $N$, the contact structure $\xi$ defines on $S$ a singular 1-foliation $S_\xi$, generated by the line field $TS\cap\xi$, that is called the {\it characteristic foliation} of $S$. Observe that generically $S$ is tangent to $\xi$ in a finite number of points that are the singularities of $S_\xi$. This one foliation is locally orientable, therefore, the index of a singular point is well defined. In the generic case the index is equal to $\pm 1$. A singular point is {\it elliptic} if the index is equal to $1$ and {\it hyperbolic} if the index is equal to $-1$.  

If $S$ is oriented or co-oriented, it is possible to induce an orientation on the characteristic foliation $S_\xi$. In this situation the singularities are endowed with a sign when we compare the orientation of $\xi$ and $TS$, that coincide in the singularity as planes. This means that positive elliptic points are sources and negative ones are sinks. For hyperbolic points the difference between positive and negative is more subtle: it is a $C^1$ rather than a topological invariant. We will always assume that $S_\xi$ is oriented.

We will distinguish two classes of contact structures: a contact structure is {\it overtwisted} if there is an embedded disc $\dd\hookrightarrow N$ whose characteristic foliation contains a limit cycle; otherwise, we will say that the contact structure is {\it tight}. The main tool for simplifying the characteristic foliation of a surface is E. Giroux's elimination lemma from \cite{giro}. Assume that we have a surface with a characteristic foliation which contains an elliptic and a hyperbolic singularities with the same sign and lying in the closure of the same leaf of the foliation. The two singularities can be eliminated via a $C^0$ small perturbation of the surface with support in a neighborhood of such a leaf. Thus we get a new surface whose characteristic foliation has two singularities less. As a non trivial corollary to the elimination lemma, in an overtwisted contact structure we can assume that there exists an embedded disc $\mathcal{D}$ whose characteristic foliation $\mathcal{D}_\xi$ has a unique elliptic singularity and the boundary $\partial\mathcal{D}$ is the limit cycle. We will call such a disc an {\it overtwisted disc}. For a visual proof of this corollary we refer to pages {\bf 28} and {\bf 29} of \cite{hofe2}.

The contact structure in the manifold $N$ can be overtwisted or tight. The proof of the theorem is based on Hofer's method to prove the existence of periodic orbits of Reeb vector fields on some closed 3-manifolds. It is divided in three cases:
\begin{itemize}
\item[(i)] when the contact structure is overtwisted;
\item[(ii)] when $\pi_2(N)\neq 0$ and the contact structure is tight;
\item[(iii)] when $N\simeq \sss^1\times \dd^2$, where $\dd^2$ is a closed disc, and the contact structure is tight.
\end{itemize}

For the proof in case (iii) we refer the reader to section {\bf 6} of \cite{etghII}. We will give the main ideas for the proof in cases (i) and (ii), see section \ref{demostracionteoremafrontera}. In the next section we will give a quick review of Hofer's method, the reader that is familiar with it can easily skip the following section.

\subsection{Pseudoholomorphic curves in H. Hofer's theorem}
\label{hofer}

The aim of this section is to sketch the proof of the following theorem in the last two situations, we follow the proof given by Hofer in \cite{hofe2}. 

\begin{teorema}[Hofer]
\label{teoremahofer2}
Let $X$ be the Reeb flow associated to a contact form $\alpha$ on a closed \mbox{3-manifold} $M$. Let $\xi$ be the contact structure defined by the kernel of $\alpha$, then $X$ has a periodic orbit in any of the following situations:
\begin{itemize}
\item $M$ is diffeomorphic to $\sss^3$;
\item $\xi$ is overtwisted;
\item $\pi_2(M)\neq 0$.
\end{itemize}
\end{teorema}

Observe that $TM=\xi\oplus X$, and  the restriction of $d\alpha$ to any plane of $\xi$ is a non degenerated 2-form. This follows from the fact that $\alpha\wedge d\alpha\neq 0$. 

We will not study the case of $\sss^3$ equipped with a tight contact structure. The proof uses the fact that all tight contact structures on $\sss^3$ are isotopic. 

The proof in the other two cases, uses pseudoholomorphic curves in a symplectisation of the manifold $M$. We are looking for periodic orbits, that we will denote by $(x, T)$ where $x:\sss^1\to M$, of the vector field $X$. Here $T$ is the minimal period of the periodic orbit. Define the functional $\Phi:C^\infty(\sss^1,M)\to \rr$ by
$$
\Phi(x)=\int_{\sss^1}x^*\alpha.
$$

\begin{proposicion}
If $x$ is a critical point of $\Phi$ and $\Phi(x)>0$, then there exists a closed integral curve $P$ of the Reeb vector field $X$ so that $x:\sss^1\to P$ is a map of positive degree. Conversely, given a closed integral curve $P$ for $X$ and a map $x:\sss^1\to P$ of positive degree, the loop $x$ is a critical point of $\Phi$ satisfying $\Phi(x)>0$.
\end{proposicion}

As we said $d\alpha$ is a non degenerated closed 2-form on the plane field $\xi$, so we can choose a compatible complex structure $J^\xi:\xi\to \xi$. The compatibility means that $d\alpha(v,J^\xi v)>0$ for every vector $v\in \xi$. The set of such complex structures is an open non empty contractible set. The manifold $M$ is now equipped with the Riemannian metric $g_{J^\xi}$ defined by
$$
g_{J^\xi}(h,k)=d\alpha(\pi(h), J^\xi\pi(k))+\alpha(h)\alpha(k),
$$
where $\pi:TM\to \xi$ is the projection along the orbits of $X$ and $h, k\in TM$.

Observe that the functional $\Phi$ and the equation $d\Phi(x)=0$ do not control the map $x$ in the $X$ direction. Such a control will be desirable for using variational methods. Formally, the \mbox{$L^2$-gradient} of the functional $\Phi$ on the loop space $C^\infty(\sss^1,M)$ associated with $d\Phi$ is the vector field $J^\xi(x)\pi\dot{x}$. The negative gradient solves the equation
$$
x=y(s) \qquad \frac{dy}{ds} =-\nabla\Phi(x),
$$
where $y:\rr\to C^\infty(\sss^1,M)$ is a smooth arc. We can define a map $v:\rr\times\sss^1\to M$, where $v(s,t)=y(s)(t)=x(t)$, that satisfies the partial differential equation
\begin{equation}\label{ecuacionxi}
\partial_sv+J^\xi(v)\pi(\partial_tv)=0.
\end{equation}
This is a first order elliptic system in the $\xi$ direction. Remark that it lacks of ellipticity in the $X$ direction. In order to control the $X$ direction, we will construct the symplectisation of $M$. Consider the non compact manifold $W=\rr\times M$ equipped with the symplectic form
$$
\omega=d(e^t\alpha)=e^t(dt\wedge\alpha+d\alpha),
$$
where $t$ is the $\rr$ coordinate. We will call $(W,\omega)$ the symplectisation of $(M,\alpha)$. Using $J^\xi$ we can define an almost complex structure $J$ on the symplectisation $W$ by
\begin{equation}\label{J}
J_{(a,h)}(b,k)=(-\alpha_h(k), J^\xi_h\pi(k)+bX_h),
\end{equation}
where $(b,k)\in T_{(a,h)}W$, and $X_h$ is the Reeb vector field on $M$ at the point $h$. Consider the Hamiltonian $H$ that is the projection from $W$ to $\rr$. The restriction of the Hamiltonian flow of $H$ to $M$ coincides with the flow of the Reeb vector field of the contact form $\alpha$.

Consider now a closed Riemann surface $(\Sigma, j)$, where $j$ is a complex structure, and take $\Gamma$ a finite set of points of $\Sigma$. A map $u:\Sigma\setminus\Gamma\to W$ is called $J$-holomorphic if
$$
du\circ j=J\circ du.
$$

\begin{lema}
If $\Gamma$ is empty the map $u$ is constant.
\end{lema}

Let $u=(a,v):\rr^+\times\sss^1\to W$, where $a:\rr^+\times\sss^1\to\rr$ is an auxiliary map. We can write equation \ref{ecuacionxi} as the next system
\begin{eqnarray}\label{va}
\pi(\partial_sv)+J^\xi(v)\pi(\partial_tv) & = & 0\\\nonumber
\alpha(\partial_tv) & = & \partial_sa\\
-\alpha(\partial_sv) & = & \partial_ta\nonumber
\end{eqnarray}
We have a first order elliptic system that controls the $X$ direction. Let us define the energy of the map $u$ as 
$$
E(u)=\sup_{f\in \Delta}\int_{\Sigma\setminus\Gamma}u^* d(f\alpha),
$$
where $\Delta=\{f:\rr\to [0,1]|f^\prime\geq 0\}$, and the one form $f\alpha$ on $\rr\times M$ is defined by
$$
(f\alpha)_{(a,h)}(b,k)=f(a)\alpha_h(k).
$$

At this point, Hofer establishes, in \cite{hofe2}, an equivalence between finding periodic orbits of $X$ on $M$ and the existence of $J$-holomorphic maps that are solutions of equation \ref{va} with finite energy. More precisely

\begin{teorema}[Hofer]{\hspace{2mm}}
Let $\Gamma\subset \Sigma$ be a finite non empty set of points. There is a finite energy non constant $J$-holomorphic map $u:\Sigma\setminus\Gamma \to W$ if and only if the Reeb vector field $X$ has a periodic orbit.
\end{teorema}

If a map $u=(a,v)$ as in the theorem exists, in the particular case where $\Sigma=\rr\times \sss^1$, we have as a consequence of the energy bound that the following limit exists
$$
\lim_{s\to\infty}\int_{\sss^1}v(s,\cdot)^*\alpha:= T\in \rr.
$$
If $T\neq 0$, then there exists a $|T|$-periodic orbit $x$ of the Reeb vector field $X$ and there exists a sequence $s_n\to\infty$, satisfying $v(s_n,t)\to x(t\;\mbox{mod}\,T)$ as $n\to\infty$, where $t\in\sss^1$ and with convergence in $C^\infty(\sss^1,M)$. The solution $x$ is then a periodic orbit of $X$ with period $T$.

For the other implication of the theorem above, assume that $X$ admits a $T$-periodic orbit $x$, that is a periodic orbit of minimal period $T$. We have to find a $J$-holomorphic map with finite energy. Define a map from the Riemann sphere minus two points to the symplectisation of $M$,
\begin{eqnarray*}
u_\pm=(a_\pm,v_\pm):\rr\times\sss^1 & \to & W\\
u_\pm(s,t) & = & (\pm Ts+c, x(\pm Tt+d))
\end{eqnarray*}
for two constants $c$ and $d$. These are $J$-holomorphic maps with zero energy, and thus solutions of the first order elliptic system \ref{va}. Clearly,
$$
\int_{\sss^1}v_\pm(s,\cdot)^*\alpha=\pm T
$$
is constant in $s\in\rr$.

We conclude, that for proving the existence of periodic orbits we need to find a finite energy non constant $J$-holomorphic map in the following two situations: when $\xi$ is an overtwisted contact structure and when $\xi$ is tight on a manifold with $\pi_2(M)\neq 0$.

\noindent{\bf Overtwisted case.}

Consider an overtwisted disc $\mathcal{D}$ in $M$, oriented in such a way that the unique elliptic singularity $e$ of $\mathcal{D}_\xi$ is positive. We can explicitly construct a one dimensional family of small $J$-holomorphic discs in $W$ with their boundaries on $\{0\}\times\mathcal{D}$ that pop out the singularity $(0,e)$. We will call such a family a Bishop family. 

\begin{teorema}{\hspace{2mm}}
\label{teoremabishop}
There is a continuous map
$$
\Psi:\dd^2\times[0,\epsilon) \to W,
$$
$\epsilon>0$, so that for each $u_t(\cdot)=\Psi(\cdot,t)$ we have that
\begin{itemize}
\item $u_t:\dd^2\to W$ is $J$-holomorphic;
\item $u_t(\partial \dd^2)\subset (\mathcal{D}\setminus \{e\})\subset \{0\}\times M$;
\item $u_t|_{\partial \dd^2}:\partial \dd^2\to (\mathcal{D}\setminus\{e\})$ has winding number $1$;
\item $\Psi|_{\dd^2\times(0,\epsilon)}$ is a smooth map;
\item $\Psi(z,0)=e$ for all $z\in \dd^2$.
\end{itemize}
\end{teorema}

It is important to notice that 
$$
u_t|_{\partial \dd^2}:\partial \dd^2\to (\mathcal{D}\setminus\{e\})
$$ 
is an embedding transversal to the characteristic foliation of $\mathcal{D}$. Following Hofer's proof (see \cite{hofe2}), we have that using the implicit function theorem we can find a {\it maximal Bishop family} 
$$
\Psi_{\max}:\dd^2\times[0,1)\to W.
$$ 
The transversality between $u_t(\partial\dd^2)$ and $\mathcal{D}_\xi$ implies that $\Psi(\partial\dd^2\times [0,1))$ cannot fill all of $\mathcal{D}$. We claim that the gradient of $\Psi_{\max}$ has to blow, that is there exist sequences $t_k\to 1$ and $z_k\to z_0\in\dd^2$ such that
$$
|\nabla\Psi_{\max}(z_k,t_k)|\to \infty.
$$
If this was not the situation, the sequence $\Psi_{\max}(\cdot, t_k)$ would converge to a $J$-holomorphic disc which will allow us to extend the maximal family $\Psi_{\max}$. This is a contradiction. Thus 
$$
|\nabla\Psi_{\max}(z_k,t_k)|\to \infty
$$
and we can assume, modulo reparameterization, that the $z_k$ are bounded away from $\partial \dd^2$. Hence the gradients are blowing up in the interior of $\dd^2$.
Let us assume that $z_k=0$ for all $k$ and that the norm of the gradient $\nabla\Psi_{\max}$ is maximal at the origin. Write 
$$
\Psi_{\max}(z,t_k)=(a_k(z),u_k(z))\in \rr\times M.
$$
Define a sequence of maps $v_k:D_k\to W$, where $D_k$ is a two dimensional disc of radius $R_k$ equal to $|\nabla\Psi_{\max}(0,t_k)|$, as
$$
v_k(z)=\left(a_k\left(\frac{z}{R_k}\right)-a_k(0),u_k\left(\frac{z}{R_k}\right)\right).
$$
The gradient of $v_k$ does not blow up. Hofer then shows that the sequence $\{v_k\}$ converge to a non constant  $J$-holomorphic finite energy plane $v:\cc=\sss^2\setminus\{\infty\}\to W$. We have constructed a finite energy $J$-holomorphic map, thus $X$ has a periodic orbit.
This finishes the proof for the overtwisted case.

\noindent{\bf The tight case where $\pi_2(M)\neq 0$.}

The sphere theorem implies that there is an embedded non contractible 2-sphere $F$ in $M$. Using again Giroux's elimination lemma, we have an embedded sphere $F$ such that $F_\xi$ has only two elliptic tangencies.

As before we can start a Bishop family of $J$-holomorphic discs at each one of the singularities. Assume that we have a uniform bound for the gradient of the two families. Under this hypothesis, we can show that the two families match up when they meet. Hence we get a continuous map $\dd\times [-1,1]\to W$ such that $\dd\times \{-1\}$ is mapped to one singularity and $\dd\times\{1\}$ to the other. That is, we get a map from the closed three dimensional ball $\dd^3$ to $W$ which induces an homeomorphism from $\sss^2=\partial \dd^3 \to F$. This implies that $F$ is contractible, which is clearly, a contradiction.

Hence we cannot have a uniform bound for the gradient of the two families, and thus, we obtain a Bishop family of $J$-holomorphic discs such that the gradients blow up. As in the overtwisted case, we can suppose that they blow up at the center of the disc and construct a \mbox{$J$-holomorphic} non constant finite energy plane which yields to the existence of a periodic orbit of $X$.

\subsection{Proof of theorem \ref{teoremafrontera}, cases (i) and (ii)}
\label{demostracionteoremafrontera}

Assume that $X|_{\partial N}$ has no periodic orbits. Note first that $\partial N$ is the union of invariant 2-tori, and the vector field $X$ is topologically conjugated to a linear vector field with irrational slope on them. Assume that $\alpha|_{N}$ is a positive contact form. Consider as in the previous section the manifold $W=\rr\times N$ equipped with the symplectic form 
$$
\omega=d(e^t\alpha)=e^t(dt\wedge\alpha+d\alpha).
$$
We will choose a complex structure $J^\xi$ on $\xi=\ker(\alpha)$, such that $d\alpha(v, J^\xi v)>0$ on $N$, for every non zero $v\in \xi$. We will use the almost complex structure $J$ defined in equation \ref{J}. The next lemma is immediate.

\begin{lema}
The boundary of $W$ is Levi flat with respect to $J$, in other words $\partial W$ is foliated by the $J$-complex surfaces $\rr\times \gamma$, where $\gamma$ is an orbit of $X$.
\end{lema}

We will use the next result by D. McDuff from \cite{mcdu} that studies the intersection between almost complex surfaces.

\begin{teorema}[McDuff]{\hspace{2mm}}
\label{teoremamcduff}
Two closed distinct $J$-holomorphic curves $C$ and $C^\prime$ in an almost complex $4$-manifold $(W,J)$ have only a finite number of intersection points. Each such a point contributes with a positive number to the algebraic intersection number $C\cdot C^\prime$.
\end{teorema}

\sdemostracion{{\bf of theorem \ref{teoremafrontera} in cases (i) and (ii)} Assume first that $\xi$ is overtwisted. We begin by completing $N$ to a closed 3-manifold $M$ and extending the contact structure $\alpha$ to $N$. Take an overtwisted disc $\mathcal{D}$ embedded in the interior of $N$. Let $W^\prime$ be the symplectisation of $M$. There exists a maximal Bishop family of $J$-holomorphic discs
$$
\Psi:\dd^2\times [0,1)\to W^\prime
$$
satisfying the conditions of theorem \ref{teoremabishop}. Observe that $\Psi(\partial\dd^2, t)\subset \mathcal{D}\subset \{0\} \times N$. We claim that $\Psi(\dd^2,t)\subset \rr\times N$. Assume that this is not the case, then one of the $u_t(\dd^2)=\Psi(\dd^2,t)$ touches the boundary of $W$ tangentially. Since $\partial W$ is foliated by $J$-holomorphic surfaces, $u_t(\dd^2)$ intersects one of the surfaces $\rr\times \gamma$, where $\gamma$ is an orbit of $X$. Since $u_t(\dd^2)$ is homotopic to a point and its boundary  through the homotopy is in the interior of $N$, the algebraic intersection number between $u_t(\dd^2)$ and $\rr\times\gamma$ is zero. Applying theorem \ref{teoremamcduff} we get a contradiction. Thus the discs in the Bishop family are inside $W$.

Recall that following the proof of Hofer we get a finite energy plane $v:\cc\to W$. Since all the $\Psi(\dd^2,t)$ are contained in the interior of $W$, so does $v(\cc)$ and thus we obtain a periodic orbit in $B$. 

The same arguments are valid when $\pi_2(N)\neq 0$. Consider a non contractible 2-sphere $F$ embedded inside $N$. Using again Giroux's elimination lemma, we have an embedded sphere $F$ such that $F_\xi$ has only two elliptic tangencies.
We can start a Bishop family of \mbox{$J$-holomorphic} discs at each one of the singularities. Using McDuff's theorem we can show that such families are contained in the symplectic manifold $\rr\times N$.

Assuming that we have a uniform bound for the gradient of the two families, we have that the two families match up when they meet. Hence we get a continuous map from $\dd\times [-1,1]\to W$ such that $\dd\times \{-1\}$ is mapped to one singularity and $\dd\times\{1\}$ to the other one. That is, we get a map from the closed three dimensional ball $\dd^3$ to $\rr\times N$ which induces an homeomorphism from $\sss^2=\partial \dd^3 \to F$. This implies that $F$ is contractible, which is clearly a contradiction.
Hence we can construct a \mbox{$J$-holomorphic} non constant finite energy plane whose image is contained in the interior of $\rr\times N$. Thus $X$ possesses a periodic orbit.
}

\subsection{Proofs of theorems \ref{teoremaperiogeodanadebil} and \ref{teoremaperiogeodvol}}

In the proofs we will also use the following proposition. To prove it we refer the reader to  the remarks about compressible and incompressible surfaces embedded in 3-manifolds, from A. Hatcher's notes \cite{hatc} (pages {\bf 11, 12}).

\begin{proposicion}\label{torosen}
Let $S$ be a finite collection of disjoint embedded tori in $\sss^3$ or a closed oriented \mbox{3-manifold} $M$ with $\pi_2(M)\neq 0$. Then there is a connected component $B$ of $\sss^3\setminus S$, respectively $M\setminus S$, such that either $\overline{B}$ is a solid torus or $\pi_2(\overline{B})\neq 0$.
\end{proposicion}

\sdemostracion{{\bf of theorem \ref{teoremaperiogeodanadebil}.} The proof is the same as section \ref{periodicageod}. We have to prove the theorem in case (III). Assume that $X$ is aperiodic, then the invariant set $A=f^{-1}(0)$ is composed by  a finite collection of invariant tori. Then by proposition \ref{torosen} there is a connected component $B$ of $M\setminus A$ whose closure is either diffeomorphic to a solid torus or $\pi_2(\overline{B})\neq 0$. Using proposition \ref{propaproximacion} we get that $X$ is the Reeb vector field of a contact form in $\overline{B}$,  and hence theorem \ref{teoremafrontera} proves the existence of a periodic orbit in $B$.
}

\sdemostracion{{\bf of theorem \ref{teoremaperiogeodvol}} Let $\omega$ be the invariant volume form, then we have a $C^\infty$ function $f$ as above that is a first integral of $X$. Again, we only need to prove the theorem in  case (III).
\begin{itemize}
\item[(III)] $f$ is equal to zero on a compact invariant set $f^{-1}(0)=A\subset M$. As we previously said $A$ is the set where the differential form $\alpha$ is closed. Observe that for a regular value $a$ of $f$, the compact set $f^{-1}(a)$ is a finite union of disjoint invariant tori. 

If zero is a regular value, $A$ is a finite union of invariant tori. Let $\epsilon$ be small enough to guarantee that the values in $[-\epsilon,\epsilon]$ are all regular. Then $f^{-1}([-\epsilon,\epsilon])$ is composed by manifolds diffeomorphic to $T\times [0,1]$ where $T$ is a two dimensional torus and the tori $T\times\{s\}$ are tangent to $X$ for every $s\in[0,1]$. Using proposition \ref{torosen}, there is a connected component of $M\setminus\{f^{-1}([-\epsilon,\epsilon])\}$ such that $\overline{B}$ is either a solid torus or a manifold with non trivial $\pi_2$. In this manifold $\overline{B}$ the vector field $X$ is tangent to the boundary and is a Reeb vector field since the restriction of $\alpha$ to $\overline{B}$ is a contact form. Thus using theorem \ref{teoremafrontera} we conclude that $X$ has a periodic orbit.

Assume now that zero is a critical value of $f$. For $\epsilon>0$ small enough let 
$$
S_\epsilon=f^{-1}(\epsilon)\cup f^{-1}(-\epsilon).
$$
Assume that $\pm \epsilon$ are regular values, then $S_\epsilon$ is a finite collection of invariant tori. 
Consider $M\setminus S_\epsilon$. By proposition \ref{torosen}, there is a connected component $B$ of $M\setminus S_\epsilon$, such that its closure is of one of the following two types:
\begin{itemize}
\item a solid torus $\sss^1\times\dd^2$, where $\dd^2$ is a two dimensional closed disc; 
\item a manifold with boundary whose second homotopy group $\pi_2$ is non trivial.
\end{itemize}
As before, denote by $\xi$ the plane field defined by the kernel of the 1-form $\alpha$. If $B\cap A=\emptyset$, we have that $\xi|_{\overline{B}}$ is a contact structure, and thus by theorem \ref{teoremafrontera} we conclude that $X$ possesses a periodic orbit in $\overline{B}$. 

We have to consider the case $B\cap A\neq \emptyset$. The plane field $\xi$ is contact in a neighborhood of $\partial\overline{B}$. We will prove that for $\epsilon$ small enough there is a closed 1-form that is transverse to $X$ in $\overline{B}$. Then by Tischler's theorem \cite{tisc} the vector field $X|_{\overline{B}}$ has a section and $\overline{B}$ is a fiber bundle over the circle. Hence $\overline{B}$ is a solid torus $\sss^1\times \dd^2$ (since a fiber bundle over the circle cannot have non trivial $\pi_2$),  and $X$ is transverse to the discs $\{\cdot\}\times \dd^2$. Then Brouwer's fixed point theorem implies that $X$ must have a periodic orbit in $\overline{B}$. 

Let $\overline{B}\subset \mathcal{B}=f^{-1}([-\epsilon,\epsilon])$. For $\epsilon$ small enough, we will construct a closed 1-form that is transverse to $X$ in $\mathcal{B}$.
We will divide the proof of the existence of the closed 1-form in two parts: first we will give explicit expressions for $X$, the forms $\alpha$ and $\iota_X\mu$ near $\mathcal{B}$, and then we will construct the closed 1-form. 

Let $0\leq \delta< \epsilon$ be small enough to guarantee that the values in the intervals $[\epsilon,\epsilon+\delta]$ and $[-\epsilon-\delta,-\epsilon]$ are all regular. Let $\mathcal{D}=f^{-1}([-\epsilon-\delta,\epsilon+\delta])$. Then $\mathcal{D}\setminus \mathcal{B}$ is foliated by tori that are tangent to $X$. Consider a connected component $D$ of $\mathcal{D}\setminus \mathcal{B}$ where $f$ is positive, and denote each invariant torus in it by $T_t$, where the $f$ equals $\epsilon+t$ on this torus and $t\in[0,\delta]$. We will do the construction in $D$ but it is analogous in the rest of $\mathcal{D}$.

On each torus there is a non singular vector field $Y$ defined by the equation $\alpha(Y)=0$ and $\iota_Y\iota_X\mu=df$.  The reason why it is non singular is that the characteristic foliation of the torus is non singular. Observe that $Y$ is tangent to each torus and is in $\xi$.

{\bf Explicit expression for $X$ in $D$}

The expression of $X$ is the same we found when proving theorem \ref{teoremaperiogeodanavol}. That is following the arguments we find that there are constants $\tau_1$ and $\tau_2$ with $\frac{\tau_1}{\tau_2}$ an irrational number, such that
$$
X=\tau_1\frac{\partial}{\partial x}+\tau_2\frac{\partial}{\partial y}.
$$

{\bf Explicit expression for $\iota_X\mu$ in $D$}

In this system of coordinates, we can write $\mu=\beta(x,y,t)dx\wedge dy\wedge dt$ for a positive function $\beta$. Then in $D$ we have that
$$
\iota_X\mu=\tau_1\beta(x,y,t)dy\wedge dt-\tau_2\beta(x,y,t)dx\wedge dt,
$$
is a closed form. Hence, $L_X\beta=0$ and since $X$ has dense orbits in each torus $T_t$, we get $\tau_1\frac{\partial \beta}{\partial x}(x,y,t)=-\tau_2\frac{\partial \beta}{\partial y}(x,y,t)$. Thus $\beta$ is just a function of $t$ and we have 
$$
\iota_X\mu=\beta(t)dt\wedge(\tau_2dx-\tau_1dy).
$$

{\bf Explicit expression for $\alpha$ in $D$}

As in the proof of theorem \ref{teoremaperiogeodanavol} we have that
$$
\alpha=\gamma(t)(\tau_2dx-\tau_1dy)+A_3(t)dt+\frac{c\tau_2dx+(1-c\tau_1)dy}{\tau_2},
$$
for some functions $\gamma$, $A_3$ and a positive constant $c$. Using the fact that $d\alpha=f\iota_X\mu$, we get that $\gamma^\prime(t)=(\epsilon+t)\beta(t)$. 

{\bf The existence of a closed 1-form in $\mathcal{B}$ transverse to $X$}

Now that we have local expressions for the forms $\alpha$ and $\iota_X\mu$, we will begin the construction of the closed 1-form. Take a $C^\infty$ function $p:[0,1]\to[0,1]$ such that $p(s)=1$ for $s<\frac{1}{3}$, for $s>\frac{2}{3}$ we set $p(s)=0$, and $p^\prime(s)\leq 0$. Define a 1-form $\tilde{\alpha}$ in $D$ as 
$$
\tilde{\alpha}=\left[\gamma(0)+p\left(\frac{t}{\delta}\right)(\gamma(t)-\gamma(0))\right](\tau_2dx-\tau_1dy)+\frac{c\tau_2dx+(1-c\tau_1)dy}{\tau_2}+A_3(t)dt,
$$
for $t\in [0,\delta]$. We can define this form in each component of $\mathcal{D}\setminus \mathcal{B}$ and extended it by $\alpha$ in $\mathcal{B}$, since $\alpha=\tilde{\alpha}$ when $t=0$. We have that
$$
d\tilde{\alpha}=\left[\frac{1}{\delta}p^\prime\left(\frac{t}{\delta}\right)(\gamma(t)-\gamma(0))+p\left(\frac{t}{\delta}\right)\gamma^\prime(t)\right]dt\wedge(\tau_2dx-\tau_1dy).
$$
In $\mathcal{D}\setminus\mathcal{B}$ the function $\gamma^\prime$ is never zero, and using the fact that $\beta(t)=\frac{\gamma^\prime(t)}{\epsilon+t}$ there is a function $h(t)$ such that $d\tilde{\alpha}=h(t)\iota_X\mu$. We have that
\begin{equation}\label{hexplicita}
h(t)=\frac{\epsilon+t}{\gamma^\prime(t)\delta}\,p^\prime\left(\frac{t}{\delta}\right)(\gamma(t)-\gamma(0))+(\epsilon+t)p\left(\frac{t}{\delta}\right),
\end{equation}
hence $h(0)=\epsilon$ and $h(\delta)=0$. Thus we get a 1-form $\tilde{\alpha}$ in $M$ such that \mbox{$d\tilde{\alpha}=h(t)\iota_X\mu$}, where $h$ equals $f$ in $\mathcal{B}$ and is equal to zero in $M\setminus \mathcal{D}$. In particular, $\tilde{\alpha}$ is closed outside $\mathcal{D}$. 

We claim that there is a positive constant $C$ independent of $\epsilon$ such that $|h|\leq C\epsilon$. In the region $D$ we have that 
$$
p\left(\frac{t}{\delta}\right)(\epsilon+t)<2\epsilon.
$$
If we choose $\delta$ small enough we can assume that $\beta(s)\leq 2\beta(t)$ for every $s\in[0,t]$ and $t\leq \delta$. Then
\begin{eqnarray*}
|\gamma(t)-\gamma(0)| & = & \Big|\int_0^t(\epsilon+s)\beta(s)ds\Big|\\
		& \leq & 2\Big|\beta(t)\left(\epsilon t+\frac{t^2}{2}\right)\Big|\\
		& \leq & 3|\beta(t)|\epsilon\delta.
\end{eqnarray*}
Putting the inequalities in equation \ref{hexplicita} we get that 
$$
|h(t)|\leq 2\epsilon\Big|p^\prime\left(\frac{t}{\delta}\right)+1\Big|\leq 2\epsilon \sup_{t} \Big|p^\prime\left(\frac{t}{\delta}\right)+1\Big|,
$$
This proves our claim.

Recall that we are looking for a 1-form in $M$ whose restriction to $\mathcal{B}$ is closed and transverse to $X$. We will now study the cohomology class of $h\iota_X\mu$ to find a 1-form different from $\alpha$ and such that its derivative is equal to $d\tilde{\alpha}$ in $\mathcal{B}$.

{\bf The cohomology class of $h\iota_X\mu$ on $M$}

Consider the exact sequence of homologies with real coefficients
$$
\cdots \to H_1(M\setminus A)\to H_1(M)\to H_1(M, M\setminus A)\to \cdots
$$
Consider a finite collection of embedded curves $\sigma_1,\sigma_2,\ldots, \sigma_n$ in $M\setminus A$ such that they form a basis for the kernel of the map $ H_1(M)\to H_1(M, M\setminus A)$. These curves are at positive distance from $A$, then for $\epsilon$ small enough we can assume that the $\sigma_i$ are at positive distance from $\mathcal{B}$.

Using the duality of Poincar\'e (see for example chapter {\bf 26} of \cite{greenb}) we have that \mbox{$H_1(M)\simeq H^2(M)$} and hence for every $i=1,2,\ldots,n$ we can find a 2-form $\omega_i$ that is the dual of $\sigma_i$ and whose support is contained in a tubular neighborhood of $\sigma_i$ contained in $M\setminus \mathcal{B}$.

\begin{lema}
For $\epsilon$ small enough there are unique real numbers $r_1,r_2,\ldots,r_n$ such that
$$
[h\iota_X\mu]=\sum_{i=1}^nr_i[\omega_i]
$$
in $H^2(M)$. Moreover, there exists a constant $C^\prime$ independent of $\epsilon$ such that $|r_i|\leq C^\prime\epsilon$ for every $i$.
\end{lema}

\demostracion{For $\epsilon$ small we can assume that $\mathcal{B}$ does not intersect the supports of the forms $\omega_i$. Denote by 
\begin{eqnarray*}
f_1:H_1(M) & \to & H_1(M,M\setminus \mathcal{B})\\
f_2:H^2(M) & \to & H^2(\mathcal{B}).
\end{eqnarray*}
Using the isomorphism given by the duality of Poincar\'e we have a map $\ker(f_1)\to \ker(f_2)$ that is injective. Recall that $h\iota_X\mu$ is exact in $\mathcal{B}$. Then to prove the existence and uniqueness of the numbers $r_i$ we need to prove that the precedent map is surjective. 

Take an element $\omega$ in the kernel of $f_2$. It can be represented by a form whose support is in $M\setminus\mathcal{B}$, then $[\omega]\in H^2_c(M\setminus\mathcal{B})$ (since it has compact support). The dual of this class under the duality of Poincar\'e is an homology class $\sigma\in H_1(M\setminus\mathcal{B})$ satisfying that for every element $S\in H_2(M\setminus\mathcal{B},\partial\mathcal{B})$
$$
\sigma\cdot S=\int_S\sigma.
$$
Using the inclusion $i:M\setminus\mathcal{B}\to M$, we get
$$
i_*\sigma\cdot S=\int_S\sigma,
$$
for all $S\in H_2(M)$. Then $i_*\sigma\in H_1(M)$ is the dual of $[\omega]\in H^2(M)$, and $f_1(i_*\sigma)=0$. Then the map is surjective.

We need to prove now that the $r_i$ are bounded. For $i=1,2,\ldots, n$ fix an oriented embedded surface $S_i$ in $M$ that intersects the $\sigma_j$. Then
$$
r_i=\int_{S_i}\sum_{j=1}^nr_j\omega_j=\int_{S_i}h\iota_X\mu.
$$
Using the bound on $h$ we get a constant $C^\prime$ that is independent of $\epsilon$ and such that $|r_i|\leq C^\prime \epsilon$.
}

The differential 2-form given by $\gamma=h\iota_X\mu-\sum_{i=1}^nr_i\omega_i$ is closed and exact in $M$. The next step is to find a primitive of $\gamma$ that is bounded by a constant multiplied by $\epsilon$.

Recall that we can define a norm on the space of $d$-forms $\Omega_d(M)$ as 
$$
\|\beta\|=\sup\{|\beta(V)|\,|\,V \, \mbox{is a unit $d$-vector}\}.
$$
The bounds above imply that $\|\gamma\|\leq C^{\prime\prime}\epsilon$ for a positive constant $C^{\prime\prime}$ independent of $\epsilon$. We need to find a primitive $\lambda$ whose norm is bounded by the norm of $\gamma$. 
The existence of such a primitive is given by combining the main result of F. Laudenbach's paper \cite{laud} and theorem {\bf 1.1} of J.-C. Sikorav's paper \cite{siko}. The first one gives a method to find a primitive and the second one a bound for it. We get,

\begin{lema}
There exists a 1-form $\lambda$ such that $d\lambda=\gamma$ and $\|\lambda\|\leq \hat{C}\|\gamma\|$, where $\hat{C}$ is a constant independent of $\epsilon$.
\end{lema}

Then, using the previous bounds we have $\|\lambda\|\leq \hat{C}C^{\prime\prime}\epsilon$. Thus the 1-form $\alpha-\lambda$ satisfies that
$$
d(\alpha-\lambda)=f\iota_X\mu-h\iota_X\mu+\sum_{i=1}^nr_i\omega_i,
$$
is equal to zero in $\mathcal{B}$, and $(\alpha-\lambda)(X)>0$ as a consequence of the bounds we found and the fact that they are independent of $\epsilon$. Then this is the 1-form we were looking for: a closed 1-form in $\mathcal{B}$ that is transverse to $X$. This finishes the proof of the theorem. 

\end{itemize}
}

Let us finish with a remark. We say that a vector field is minimal if all its orbits are dense in the ambient manifold. The still open W. H. Gottschalk conjecture asserts that there are no minimal vector fields on $\sss^3$. Observe that a geodesible vector field on $\sss^3$ cannot be minimal, in fact the only minimal geodesible vector fields on closed 3-manifolds are the suspensions of minimal diffeomorphisms of a two dimensional torus. To prove this claim consider a minimal geodesible vector field on a closed 3-manifold $M$. Then the invariant set $A=f^{-1}(0)$ must be equal to $M$ or empty. In the latter case the vector field is a Reeb vector field of a contact structure, then it cannot be minimal since it possesses a periodic orbit. If $A=M$, the vector field admits a global section that must be a torus since $X$ is aperiodic. Then it is the suspension of a minimal diffeomorphism of a two dimensional torus and $M$ is a torus bundle over the circle.

The difficulty in extending theorem \ref{teoremaperiogeodvol} to geodesible vector fields that are not volume preserving is dealing with the set $A$. Thus the question whether a geodesible vector field on $\sss^3$ has a periodic orbit is still open. The more general question is to find necessary and sufficient conditions for the existence of periodic orbits of non singular vector fields on closed 3-manifolds.

\small{

}

\end{document}